# Nice q-analogs of orthogonal polynomials with nice moments: Some simple examples.

Johann Cigler

**Abstract**

In this note I collect some typical examples of orthogonal polynomials with simple moments where both moments and orthogonal polynomials have nice $q$ – analogs.

## 1. Introduction

Let $(a(n))_{n\geq 0}$ be a sequence of real numbers with $a(0)=1$. For the computation of the Hankel determinants $d_n = \det(a(i+j))_{i,j=0}^{n-1}$ it is often useful to consider the linear functional $L$ on $\mathbb{R}[x]$ defined by $L(x^n) = a(n)$. If $d_n \neq 0$ for all $n \in \mathbb{N}$ there exists a sequence of monic polynomials $p_n(x)$ which are orthogonal with respect to $L$, i.e. satisfy $L(p_k(x)p_n(x)) = 0$ for $k \neq n$ and $L(p_n(x)^2) \neq 0$. In some cases there exist nice $q$ – analogs $a(n;q)$ of $a(n)$ where the corresponding orthogonal polynomials $p_n(x;q)$ are nice $q$ – analogs of $p_n(x)$ too. In this note I collect some simple examples.

## 2. Some background material

Let $(a(n))_{n\geq 0}$ be a sequence of real numbers with $a(0)=1$ such that no Hankel determinant vanishes, and let $L$ be the linear functional on $\mathbb{R}[x]$ with moments $L(x^n) = a(n)$.

Then the following situation is well known (cf. [6]). There exist monic polynomials $p_n(x)$ which are orthogonal with respect to $L$. They are given by

(1) $$p_n(x) = \frac{1}{d_n} \det \begin{pmatrix} a(0) & a(1) & \ldots & a(n-1) & 1 \\ a(1) & a(2) & \ldots & a(n) & x \\ a(2) & a(3) & \ldots & a(n+1) & x^2 \\ \vdots & \vdots & \vdots & \vdots & \vdots \\ a(n) & a(n+1) & \ldots & a(2n-1) & x^n \end{pmatrix}.$$

To see this, observe that for $k < n$

(2) $$L(x^k p_n(x)) = \frac{1}{d_n} \det \begin{pmatrix} a(0) & a(1) & \ldots & a(n-1) & L(x^k) \\ a(1) & a(2) & \ldots & a(n) & L(x^{k+1}) \\ a(2) & a(3) & \ldots & a(n+1) & L(x^{k+2}) \\ \vdots & \vdots & \vdots & \vdots & \vdots \\ a(n) & a(n+1) & \ldots & a(2n-1) & L(x^{k+n}) \end{pmatrix} = 0,$$

because two columns are equal and that $L(x^n p_n(x)) = d_{n+1} \neq 0$.

By a theorem of Favard there exist numbers $s_n$ and $t_n$ such that

(3) $$p_n(x) = (x - s_{n-1})p_{n-1}(x) - t_{n-2}p_{n-2}(x).$$

The coefficients $a(n,k)$ in the expansion

(4) $$x^n = \sum_{k=0}^{n} a(n,k) p_k(x)$$

are uniquely determined by

(5) $$a(n,k) = a(n-1,k-1) + s_k a(n-1,k) + t_k a(n-1,k+1)$$

with $a(n,k) = 0$ for $k < 0$ and $a(0,k) = [k = 0]$.

Of course, we have $a(n,0) = L(x^n) = a(n)$.

The Hankel determinants depend only on $t_n$ (cf. e.g. [2])

(6) $$d_n = \det\left(a(i+j)\right)_{i,j=0}^{n-1} = \prod_{i=1}^{n-1}\prod_{j=0}^{i-1} t_j.$$

In this note I collect some typical examples of "naturally" occurring sequences $(a(n))$ with "nice" $q$–analogs $(a(n;q))$ whose orthogonal polynomials $p_n(x;q)$ are "nice" $q$–analogs of $p_n(x)$ too, where moreover $s_n$ and $t_n$ can easily be guessed by computing $p_n(x)$ for some small values of $n$. All results are either well known or special cases of general theorems and can also easily be verified. Therefore, I can leave the proofs to the reader.

For $q$–analogs we use the standard notations $[n] = [n]_q = 1 + q + \cdots + q^{n-1} = \dfrac{1 - q^n}{1 - q}$,

$[n]! = [n]_q! = \prod_{j=1}^{n} [j]_q$, $\begin{bmatrix} n \\ k \end{bmatrix}_q = \dfrac{[n]!}{[k]![n-k]!}$ for $k \in \{0, 1, \cdots, n\}$ and

$(x;q)_n = (1-x)(1-qx)\cdots(1-q^{n-1}x)$. Here $q$ can be a real number or an indeterminate.

For given sequences $(a(n))$ we shall also consider the aerated sequences $(A(n))$ with $A(2n) = a(n)$ and $A(2n+1) = 0$. Their orthogonal polynomials $P_n(x)$ satisfy

(7) $$P_n(x) = xP_{n-1}(x) - T_{n-2}P_{n-2}(x)$$

for some sequence $(T_n)_{n \geq 0}$.

The numbers $s_n$ and $t_n$ for the original sequence $(a(n))$ are then given by

(8) $$\begin{aligned} s_n &= T_{2n-1} + T_{2n}, \\ t_n &= T_{2n}T_{2n+1} \end{aligned}$$

and their orthogonal polynomials by

(9) $$p_n(x) = P_{2n}(\sqrt{x}).$$

**Remark**

Let us note a degenerate case: For the sequence $a(n) = 1$ and the linear functional $L$ defined by $L(x^n) = 1$ the polynomials $p_n(x) = (x-1)^n$ satisfy $L(p_m(x)p_n(x)) = 0$ for all $(m,n) \neq (0,0)$. For the $q$-analog $a(n,q) = q^{\binom{n}{2}}$ of $a(n)$ the uniquely defined orthogonal polynomials are the nice $q$-analogs of $(x-1)^n$

(10) $$p_n(x,q) = \sum_{j=0}^{n} (-1)^j q^{(n-1)j} \begin{bmatrix} n \\ j \end{bmatrix}_q x^{n-j}$$

and the numbers $a(n,k) = q^{\binom{n}{2} - \binom{k}{2}} \begin{bmatrix} n \\ k \end{bmatrix}_q$ are nice $q$-analogs of the binomial coefficients.

To show the orthogonality let $F$ be the linear functional defined by $F(x^n) = q^{\binom{n}{2}}$. Then

$$F(x^m p_n(x)) = \sum_{j=0}^{n} (-1)^j q^{(n-1)j} \begin{bmatrix} n \\ j \end{bmatrix}_q F(x^{m+n-j}) = \sum_{j=0}^{n} (-1)^j q^{(n-1)j} \begin{bmatrix} n \\ j \end{bmatrix}_q q^{\binom{m+n-j}{2}}$$

$$= q^{\binom{n}{2} + \binom{m}{2} + mn} \sum_{j=0}^{n} (-1)^j q^{\binom{j}{2}} q^{-jm} \begin{bmatrix} n \\ j \end{bmatrix}_q = q^{\binom{n}{2} + \binom{m}{2}} (q^m - 1)(q^m - q) \cdots (q^m - q^{n-1}).$$

using Rothe's Theorem $\sum_{j=0}^{n} (-1)^j q^{\binom{j}{2}} \begin{bmatrix} n \\ j \end{bmatrix}_q x^j = (1-x)(1-qx)\cdots(1-q^{n-1}x).$

The right-hand side vanishes for $n > m$. For $n = m$ we get $q^{3\binom{n}{2}}(q-1)^n[n]!$

## 3. Factorials and Laguerre polynomials.

Let $a(n) = \dfrac{(n+m)!}{m!}$ for $m \in \mathbb{N}$.

For the aerated sequence $(A(n))$ we get

(11) $$\begin{aligned} T_{2n} &= n+1+m, \\ T_{2n+1} &= n+1 \end{aligned}$$

and

(12)
$$P_{2n}(x) = \sum_{j=0}^{n}(-1)^j \binom{n}{j}\binom{n+m}{j} j! x^{2n-2j},$$
$$P_{2n+1}(x) = \sum_{j=0}^{n}(-1)^j \binom{n}{j}\binom{n+m+1}{j} j! x^{2n+1-2j}.$$

This implies

(13)
$$t_n = (n+1)(n+1+m),$$
$$s_n = 2n+1+m.$$

The orthogonal polynomials are the monic Laguerre polynomials

(14)
$$p_n(x) = L_n^{(m)}(x) = \sum_{j=0}^{n}(-1)^{n-j}\binom{n+m}{n-j}\frac{n!}{j!}x^j = \sum_{j=0}^{n}(-1)^j\binom{n}{j}\frac{a(n)x^{n-j}}{a(n-j)}$$

The numbers $a(n,k)$ are

(15)
$$a(n,k) = \binom{n}{k}\frac{a(n)}{a(k)}.$$

A natural $q$–analog of $a(n)$ is $a(n;q) = \dfrac{[n+m]_q!}{[m]_q!}.$

Considering the linear functional $L$ on $\mathbb{R}(q)[x]$ defined by $L(x^n) = a(n;q)$ we get

(16)
$$T_{2n} = q^n[n+1+m],$$
$$T_{2n+1} = q^{n+m+1}[n+1]$$

and thus

(17)
$$t_n = q^{2n+m+1}[n+1]_q[n+1+m]_q,$$
$$s_n = q^n\left([n+m+1] + q^m[n]\right).$$

The orthogonal polynomials are the monic $q$–Laguerre polynomials

(18)
$$p_n(x;q) = L_n^{(m)}(x;q) = \sum_{k=0}^{n}(-1)^{n-k}q^{\binom{n-k}{2}}\frac{[n]!}{[k]!}\begin{bmatrix}n+m\\n-k\end{bmatrix}x^k$$
$$= \sum_{k=0}^{n}(-1)^k q^{\binom{k}{2}}\begin{bmatrix}n\\k\end{bmatrix}_q \frac{a(n;q)x^{n-k}}{a(n-k;q)}$$

and the $a(n,k)$ are given by

(19) $$a(n,k;q) = \frac{[n+m]!}{[m+k]!}\begin{bmatrix}n\\k\end{bmatrix}_q = \begin{bmatrix}n\\k\end{bmatrix}_q \frac{a(n;q)}{a(k;q)}.$$

**Remark**

Since

(20) $$a\left(n;\frac{1}{q}\right) = \frac{[n+m]!}{[m]!}\frac{1}{q^{\binom{n}{2}+mn}} = \frac{a(n;q)}{q^{\binom{n}{2}+mn}}$$

and

(21) $$p_n\left(x;\frac{1}{q}\right) = a\left(n;\frac{1}{q}\right)\sum_{k=0}^{n}(-1)^k q^{-kn+\binom{k+1}{2}}\begin{bmatrix}n\\k\end{bmatrix}_q \frac{x^{n-k}}{a\left(n-k;\frac{1}{q}\right)}$$

we see that $\dfrac{a(n;q)}{q^{\binom{n}{2}}}$ also is a nice $q$-analog with nice orthogonal polynomials. But it seems that $k=0$ and $k=-1$ are the only values of $k$ such that $a(n;q)q^{k\binom{n}{2}}$ has nice orthogonal polynomials.

### 4. A generalization to multifactorials

Define the multifactorial $mf(n,r) = n!^r$ by $mf(n,r) = n*mf(n-r,r)$ with $mf(n,r) = 1$ for $n \leq 1$ and consider

(22) $$a(n) = \frac{mf(rn+m,r)}{mf(m,r)} = \frac{(rn+m)!^r}{(m)!^r}.$$

Then we get

(23) $$\begin{aligned}T_{2n} &= r(n+1)+m,\\ T_{2n+1} &= r(n+1).\end{aligned}$$

(24) $$\begin{aligned}t_n &= r(n+1)(rn+r+m),\\ s_n &= (2n+1)r+m.\end{aligned}$$

(25) $$p_n(x) = \sum_{k=0}^{n}(-1)^k \binom{n}{k}\frac{a(n)}{a(n-k)}x^{n-k}$$

and

(26) $$a(n,k) = \binom{n}{k}\frac{a(n)}{a(k)}.$$

For the $q$ – analog we define the multifactorial $mf(n,r,q) = [n]!^r$ by
$mf(n,r,q) = [n] \cdot mf(n-r,r,q)$ with $mf(n,r,q) = 1$ for $n \leq 1$.

We consider $a(n;q) = \dfrac{[rn+m]!^r}{[m]!^r}$. Here we get

(27) $$T_{2n} = q^{rn}[r(n+1)+m],$$
$$T_{2n+1} = q^{r(n+1)+m}[r(n+1)].$$

(28) $$t_n = [(n+1)r][r(n+1)+m]q^{r(2n+1)+m},$$
$$s_n = q^{rn}\left([r(n+1)+m]+q^m[rn]\right)$$

and

(29) $$p_{n,rm,r}(x;q) = \sum_{k=0}^{n}(-1)^k q^{r\binom{k}{2}}\begin{bmatrix}n\\k\end{bmatrix}_{q^r}\begin{bmatrix}n+m\\k\end{bmatrix}_{q^r}[rk]!^r\, x^{n-k},$$

$$p_{n,m,r}(x;q) = \sum_{k=0}^{n}(-1)^k q^{r\binom{k}{2}}\begin{bmatrix}n\\k\end{bmatrix}_{q^r}\frac{a(n,q)}{a(n-k,q)}x^{n-k}.$$

**Remark**

For $r=2$ and $a(n) = (2n-1)!!$ the formulae can be simplified due to the identity
$\begin{bmatrix}2n\\2k\end{bmatrix}_q[2k-1]!! = \begin{bmatrix}n\\k\end{bmatrix}_{q^2}\dfrac{[2n-1]!!}{[2n-2k-1]!!}$. Then we get the Hermite polynomials

(30) $$P_n(x) = H_n(x) = \sum_{k=0}^{\lfloor n/2 \rfloor}(-1)^k\binom{n}{2k}(2k-1)!!\, x^{n-2k},$$

$$p_n(x) = h_n(x) = \sum_{k=0}^{n}(-1)^k\binom{2n}{2k}(2k-1)!!\, x^{n-k}$$

and

(31) $$a(n,k) = \frac{(2n-1)!!}{(2k-1)!!}\binom{n}{k}.$$

For $a(n;q) = [2n-1]!!$ we get the discrete $q$ – Hermite polynomials

(32) $$p_n(x;q) = h_n(x,q) = \sum_{k=0}^{n}(-1)^k q^{2\binom{k}{2}}\begin{bmatrix}2n\\2k\end{bmatrix}_q[2k-1]!!\, x^{n-k}$$

and

(33) $$a(n,k;q) = \begin{bmatrix} n \\ k \end{bmatrix}_{q^2} \frac{[2n-1]!!}{[2k-1]!!}.$$

## 5. Catalan numbers and central binomial coefficients

It is well known that the Catalan numbers $C_n = \frac{1}{n+1}\binom{2n}{n}$ are the even moments $C_n = L(x^{2n})$ of the Fibonacci polynomials $F_n(x) = \sum_{k=0}^{\lfloor n/2 \rfloor} (-1)^k \binom{n-k}{k} x^{n-2k}$ which satisfy the recursion $F_n(x) = xF_{n-1}(x) - F_{n-2}(x)$ with initial values $F_0(x) = 1$ and $F_1(x) = x$ and the central binomial coefficients $\binom{2n}{n}$ are the moments of the Lucas polynomials

$$L_n(x) = \sum_{k=0}^{\lfloor n/2 \rfloor} (-1)^k \frac{n}{n-k} \binom{n-k}{k} x^{n-2k} \text{ with } L_0(x) = 1.$$

There are many known $q$-analogs of the Catalan numbers (cf.[5]). The most natural one would be $\frac{1}{[n+1]} \begin{bmatrix} 2n \\ n \end{bmatrix}$. Unfortunately, the corresponding orthogonal polynomials (which of course are $q$-analogs of the Fibonacci polynomials) don't have nice coefficients.

The situation looks somewhat different if we replace the Fibonacci and Lucas polynomials by the corresponding monic Chebyshev polynomials. As is well known the Chebyshev polynomials of the second kind $U_n(x)$ satisfy $U_n(x) = 2xU_{n-1}(x) - U_{n-2}(x)$ with initial values $U_0(x) = 1$ and $U_1(x) = 2x$ and are therefore related to the Fibonacci polynomials by $U_n(x) = F_n(2x)$. The monic versions $u_n(x)$ satisfy

(34) $$u_n(x) = \frac{U_n(x)}{2^n} = \sum_{k=0}^{\lfloor n/2 \rfloor} (-1)^k \binom{n-k}{k} \frac{1}{4^k} x^{n-2k}.$$

Their moments $L(x^{2n}) = \frac{C_n}{4^n}$ can be written as $\frac{C_n}{4^n} = 2\frac{(2n-1)!!}{(2n+2)!!}.$

With this formulation a natural $q$-analog would be

(35) $$[2]\frac{[2n-1]!!}{[2n+2]!!} = \frac{1}{[n+1]_{q^2}} \begin{bmatrix} 2n \\ n \end{bmatrix}_{q^2} \frac{1}{\prod_{j=1}^{2n}(1+q^j)}.$$

The right-hand side is a special case of the $q$-Catalan numbers which have been proposed by George Andrews [1] from the point of view of hypergeometric series. Surprisingly, it turns out that these $q$-analogs do have orthogonal polynomials with nice coefficients.

For the aerated sequence we get

(36) $$T_n = \frac{q^n}{(1+q^{n+1})(1+q^{n+2})}.$$

The corresponding monic orthogonal polynomials $u_n(x,q)$ are

(37) $$P_n(x) = u_n(x,q) = \sum_{j=0}^{\lfloor \frac{n}{2} \rfloor} (-1)^j q^{2\binom{j}{2}} \begin{bmatrix} n-j \\ j \end{bmatrix}_{q^2} \frac{1}{(-q^{n-2j+1},q)_{2j}} x^{n-2j}.$$

Note that $u_n(x,q)$ satisfies $u_n(x,q) = xu_{n-1}(x,q) - \frac{q^{n-2}}{(1+q^{n-1})(1+q^n)} u_{n-2}(x,q)$ with initial values $u_0(x,q) = 1$ and $u_1(x,q) = x$. If we set $U_n(x,q) = (1+q)(1+q^2)\cdots(1+q^n)u_n(x,q)$ we get

$U_n(x,q) = (1+q^n)U_{n-1}(x,q) - q^{n-2}U_{n-2}(x,q)$ with $U_0(x,q) = 1$ and $U_1(x,q) = (1+q)x$.

Thus $U_n(x,q)$ is a nice $q-$analog of the Chebyshev polynomials of thé second kind. These polynomials have been introduced in [2] with other methods.

The Chebyshev polynomials of the first kind $T_n(x)$ satisfy $T_n(x) = 2xT_{n-1}(x) - T_{n-2}(x)$ with initial values $T_0(x) = 1$ and $T_1(x) = x$ and are therefore related to the Lucas polynomials by $2T_n(x) = L_n(2x)$. The monic versions $t_n(x)$ satisfy

(38) $$t_n(x) = \frac{T_n(x)}{2^{n-1}} = \sum_{k=0}^{\lfloor \frac{n}{2} \rfloor} (-1)^k \frac{n}{n-k} \binom{n-k}{k} \frac{1}{4^k} x^{n-2k}.$$

Their even moments $L(x^{2n}) = \frac{1}{4^n}\binom{2n}{n}$ can we written as $\frac{1}{4^n}\binom{2n}{n} = \frac{(2n-1)!!}{(2n)!!}$.

With this formulation a natural $q-$analog would be

(39) $$\frac{[2n-1]!!}{[2n]!!} = \begin{bmatrix} 2n \\ n \end{bmatrix}_q \frac{1}{(-q;q)_n(-q^{n+1};q)_n}.$$

For the aerated sequence we get

(40)
$$T_0 = \frac{1}{1+q},$$
$$T_n = \frac{q^n}{(1+q^n)(1+q^{n+1})}.$$

The monic orthogonal polynomials are then given by

(41)
$$t_n(x,q) = \sum_{k=0}^{\lfloor \frac{n}{2} \rfloor} (-1)^k q^{2\binom{k}{2}} \frac{[n]}{[n-k]} \begin{bmatrix} n-k \\ k \end{bmatrix}_q \frac{1}{(-q;q)_k (-q^{n-k};q)_k} x^{n-2k}$$

and satisfy $t_n(x,q) = x t_{n-1}(x,q) - \frac{q^{n-2}}{(1+q^{n-2})(1+q^{n-1})} t_{n-1}(x,q)$ for $n \geq 3$ and

$t_2(x,q) = x t_1(x,q) - \frac{1}{1+q} t_0(x,q)$, with initial values $t_0(x,q) = 1$ and $t_1(x,q) = x$.

If we set $T_n(x,q) = (-q;q)_{n-1} t_n(x,q)$ for $n \geq 1$ and $T_0(x,q) = 1$ we get

$T_n(x,q) = (1+q^{n-1}) x T_{n-1}(x,q) - q^{n-2} T_{n-2}(x,q)$ with $T_0(x,q) = 1$ and $T_1(x,q) = x$.

Thus $T_n(x,q)$ is a nice $q$ – analog of the Chebyshev polynomials of the first kind.

**Remark**

It should be noted that for the (non-orthogonal) $q$ – Fibonacci polynomials

(42)
$$F_n(x,q) = \sum_{k=0}^{\lfloor \frac{n}{2} \rfloor} (-1)^k q^{\binom{k}{2}} \begin{bmatrix} n-k \\ k \end{bmatrix} x^{n-2k}$$

the even moments of the linear functional defined by $L(F_n(x,q)) = [n = 0]$ are the $q$ – Catalan numbers $C_n(q) = \frac{1}{[n+1]} \begin{bmatrix} 2n \\ n \end{bmatrix}$.

Similarly for the (non-orthogonal) $q$ – Lucas polynomials

(43)
$$L_n(x,q) = \sum_{k=0}^{\lfloor \frac{n}{2} \rfloor} (-1)^k q^{\binom{k}{2}} \begin{bmatrix} n-k \\ k \end{bmatrix} \frac{[n]}{[n-k]} x^{n-2k}$$

the even moments of the linear functional defined by $L(L_n(x,q)) = [n = 0]$ are the $q$ – central binomial coefficients $\begin{bmatrix} 2n \\ n \end{bmatrix}$.